\documentclass{amsart}
\newtheorem{theorem}{Theorem}
\newtheorem*{theorema}{Theorem A}
\begin{document}
\title{Covering dimension and nonlinear equations}
\author{Biagio Ricceri}
\address{Department of Mathematics\\
University of Catania\\
Viale A. Doria 6\\
95125 Catania, Italy}
\maketitle

For a set $S$ in a Banach space, we denote by $\dim(S)$ its covering 
dimension \cite[p.~42]{1}. Recall that, when $S$ is a convex set, the 
covering dimension of $S$ coincides with the algebraic dimension of $S$, 
this latter being understood as $\infty$ if it is not finite 
\cite[p.~57]{1}. Also, $\overline{S}$ and $\operatorname{conv}(S)$ will 
denote the closure and the convex hull of $S$, respectively.

In [3], we proved what follows.

\begin{theorema}[{\cite[Theorem 1]{3}}]
Let $X, Y$ be two Banach spaces, $\Phi\colon X\to Y$ a continuous, 
linear, surjective operator, and $\Psi\colon X\to Y$ a continuous 
operator with relatively compact range. Then, one has
$$\dim(\{x\in X : \Phi(x)=\Psi(x)\}) 
\geq \dim(\Phi^{-1}(0)).$$
\end{theorema}
In the present paper, we improve Theorem A by establishing the
following result.

\begin{theorem}
Let $X, Y$ be two Banach spaces, $\Phi\colon X\to Y$ a continuous, 
linear, surjective operator, and $\Psi\colon X\to Y$ a completely 
continuous operator with bounded range. Then, one has
$$\dim(\{x\in X:\Phi(x)=\Psi(x)\})
\geq \dim(\Phi^{-1}(0)).$$
\end{theorem}

\begin{proof}
First, assume that $\Phi$ is not injective. For each  $x\in X$, $y\in Y$, 
$r>0$, we denote by $B_{X}(x,r)$ (resp. $B_{Y}(y,r)$) the closed ball in 
$X$ (resp. $Y$) of radius $r$ centered at $x$ (resp. $y$). By the open
mapping theorem, there is $\delta>0$ such that
$$B_{Y}(0,\delta)\subseteq \Phi(B_{X}(0,1)).$$
Since $\Psi(X)$ is bounded, there is $\rho>0$ such that
$$\bar{\Psi(X)} \subseteq B_{Y}(0,\rho).$$
Consequently, one has
$${\overline {\Psi(X)}}\subseteq
 \Phi\left (B_{X}\left (0,{{\rho}\over {\delta}}\right )\right ).$$
Now, fix any bounded open convex set $A$ in $X$ such that 
$$B_{X}\left(0,{{\rho}\over {\delta}}\right)\subseteq A.$$
Put
$$K={\overline {\Psi(A)}}.$$
Since $\Psi$ is completely continuous, $K$ is compact.
Fix any positive integer $n$ such that $n\leq \dim(\Phi^{-1}(0))$.
Also, fix $z\in K$.
Thus, $\Phi^{-1}(z)\cap A$ is a convex set of dimension at least $n$.
Choose $n+1$ affinely independent points $u_{z,1},\dots,u_{z,n+1}$ in
$\Phi^{-1}(z)\cap A$. By the open mapping theorem again, the operator 
$\Phi$ is open, and so, successively, the multifunctions 
$y\to \Phi^{-1}(y)$, $y\to \Phi^{-1}(y)\cap A$, and 
$y\to {\overline {\Phi^{-1}(y)\cap A}}$ are lower semicontinuous. 
Then, applying the classical Michael theorem \cite[P.~98]{2} to the 
restriction to $K$ of the latter multifunction, we get $n+1$ continuous
functions $f_{z,1},\dots,f_{z,n+1}$, from $K$ into ${\overline A}$, such 
that, for all $y\in K, i=1,\dots,n+1$, one has
$$\Phi(f_{z,i}(y))=y$$
and
$$f_{z,i}(z)=u_{z,i}.$$
Now, for each $i=1,\dots,n+1$, fix a neighbourhood $U_{z,i}$ of $u_{z,i}$ 
in $A$ in such a way that, for any choice of $w_{i}$ in $U_{z,i}$, the 
points $w_{1},\dots,w_{n+1}$ are affinely independent.
Now, put
$$V_{z}=\bigcap_{i=1}^{n+1}f_{z,i}^{-1}(U_{z,i}).$$
Thus, $V_{z}$ is a neighbourhood of $z$ in $K$. Since $K$ is compact, 
there are finitely many $z_{1},\dots,z_{p}\in K$ such that 
$K=\cup_{j=1}^{p}V_{z_{j}}$.
For each $y\in K$, put
$$F(y)
= \operatorname{conv}(\{f_{z_{j},i}(y) : j=1,\dots,p, i=1,\dots,n+1\}).$$
Observe that, for some $j$, one has $y\in V_{z_{j}}$, and so
$f_{z_{j},i}(y)\in U_{z_{j},i}$ for all $i=1,\dots,n+1$. Hence, $F(y)$ is 
a compact convex subset of $\Phi^{-1}(y)\cap {\overline A}$, with 
$\dim(F(y))\geq n$.
Observe also that the multifunction is $F$ is continuous 
(\cite[p.~86 and p.~89]{2}) and that the set $F(K)$ is compact 
(\cite[p.~90]{2}). Put 
$$C=\operatorname{conv}(F(K)).$$
Furthermore, note that, by continuity, one has 
$\Psi({\overline A}) \subseteq K$. Finally, consider the multifunction 
$G\colon C\to 2^{C}$ defined by putting
$$G(x)=F(\Psi(x))$$
for all $x\in C$. Hence, $G$ is a continuous multifunction, from the 
compact convex set $C$ into itself, whose values are compact convex
sets of dimension at least $n$. Consequently, by the result of \cite{4}, 
one has
$$\dim(\{x\in C : x\in G(x)\})\geq n.$$
But, since
$$\{x\in C : x\in F(\Psi(x))\}
\subseteq \{x\in C : x\in \Phi^{-1}(\Psi(x))\}$$
the conclusion follows (\cite[p.~220]{1}).
Finally, if $\Phi$ is injective, the conclusion means simply that the set 
$\{x\in X : \Phi(x)=\Psi(x)\}$ is non-empty, and this is got readily 
proceeding as before.
\end{proof}

In \cite{3}, we indicated some examples of application of Theorem A. 
We now point out an application of Theorem 1 which cannot be obtained 
from Theorem A. For a Banach space $E$, we denote by $\mathcal{L}(E)$ the
space of all continuous linear operators from $E$ into $E$, with the 
usual norm. Also, $I$ will denote a (non-degenerate) compact real 
interval.

\begin{theorem}
Let $E$ be an infinite-dimensional Banach space, $A\colon I\to 
\mathcal{L}(E)$ a continuous function and $f\colon I\times E\to E$ a 
uniformly continuous function with relatively compact range.
Then, one has
$$\dim(\{u\in C^{1}(I,E) : u'(t)=A(t)(u(t))+f(t,u(t)), 
\forall t\in I\})
= \infty.$$
\end{theorem}

\begin{proof}
Take $X=C^{1}(I,E)$, $Y=C^{0}(I,E)$ and 
$\Phi(u)=u'(\cdot)-A(\cdot)(u(\cdot))$ for all $u\in X$. 
So, by a classical result, $\Phi$ is a continuous linear operator from 
$X$ onto $Y$ such that $\dim(\Phi^{-1}(0))=\infty$. Next, put 
$\Psi(u)=f(\cdot,u(\cdot))$ for all $u\in X$. So, $\Psi$ is an operator 
from $X$ into $Y$ with bounded range. From our assumptions, thanks to the 
Ascoli-Arzel\`a theorem, it also follows that $\Psi$ is completely 
continuous. Then, the conclusion follows directly from Theorem 1.
\end{proof}

Analogously, one gets from Theorem 1 the following

\begin{theorem}
Let $A\colon I\to \mathcal{L}({\bf R}^n)$ be a continuous function and 
$f\colon I\times {\bf R}^{n}\to {\bf R}^n$ a continuous and bounded 
function. Then, one has
$$\dim(\{u\in C^{1}(I,{\bf R}^n) : u'(t)=A(t)(u(t))+f(t,u(t),
\forall t\in I\})\geq n.$$
\end{theorem}

\begin{theorem}
Let $a_{1},\dots,a_{k}$ be $k$ continuous real functions on $I$. 
Further, let $f\colon I\times {\bf R}^{k}\to {\bf R}$ be a continuous and 
bounded function.
Then, one has
$$\textstyle
\dim
\left(\left\{
u\in C^{k}(I) : 
u^{(k)}(t)
+ \sum_{i=1}^{k}a_{i}(t)u^{(k-i)}(t)
= f(t,u(t),u'(t),\dots,u^{(k-1)}(t)),
\forall t\in I \right\}\right)
\geq k.$$
\end{theorem}

\end{document}